\documentclass[12pt]{article}
\usepackage{amssymb,latexsym,start2e,graphs2,pstcol,pst-plot}

\newcommand{\Mon}{\mathop{\rm Mon}\nolimits}

\begin{document}
\pagestyle{plain}

\title{Monomial Bases for Broken Circuit Complexes
}
\author{
Jason I. Brown\\[-5pt]
\small Department of Mathematics and Statistics\\[-5pt] 
\small Dalhousie University\\[-5pt]
\small Halifax, NS B3H 3J5\\[-5pt] 
\small CANADA\\[-5pt]
\small \texttt{brown@mscs.dal.ca}\\[10pt]
Bruce E. Sagan\\[-5pt]
\small Department of Mathematics\\[-5pt] 
\small Michigan State University\\[-5pt]
\small East Lansing, MI 48824-1027\\[-5pt] 
\small USA\\[-5pt]
\small \texttt{sagan@math.msu.edu}\\[10pt]
\it This paper is dedicated to Michel Las Vergnas\\
\it on the occasion of his 65th birthday.
}

\date{\today\\[10pt]
      \begin{flushleft}
      \small Key Words: acyclic orientation, cocircuit, $f$-vector, $h$-vector,
      homogeneous system of parameters, monomial basis, no broken
      circuit, order ideal, Stanley-Reisner ring, theta graph, upper bound\\[5pt]
      \small AMS subject classification (2000):
      Primary 05E99;
      Secondary 05C99, 13A02.
      \end{flushleft}
     }
\maketitle

\begin{abstract}
Let $F$ be a field and let $G$ be a finite graph with a total ordering
on its edge set.   Richard Stanley noted that the Stanley-Reisner
ring $F(G)$ of the broken circuit complex of $G$ is 
Cohen-Macaulay.  Jason Brown gave an explicit description of a
homogeneous system of parameters for $F(G)$ in terms of fundamental
cocircuits in $G$.  So $F(G)$ modulo this hsop is a finite
dimensional vector space.  We conjecture an explicit monomial basis
for this vector  
space in terms of the circuits of $G$ and prove that the conjecture
is true for two infinite families of graphs.  We also explore an
application of these ideas to bounding the number of acyclic
orientations of $G$ from above. 
\end{abstract}

\section{Simplicial complexes and chromatic polynomials}

Let $E$ be a finite set and let $\De$ be an {\it abstract simplicial
complex on $E$}, i.e., a nonempty family of subsets of $E$ such that $S\in\De$
and $T\sbe S$ implies $T\in\De$.  The elements $S$ of $\De$ are called
{\it faces\/}.  We will assume henceforth that $\De$ is 
{\it pure of rank $r$\/} which means that all maximal faces $S$ have
$|S|=r$  where the absolute value sign denotes cardinality.
Let $f_i=f_i(\De)$ be the number of $S\in\De$ with $|S|=i$.  
Then $\De$ has {\it $f$-vector\/}
$$
\bff=\bff(\De)=(f_0,f_1,\ldots,f_r)
$$
as well as {\it $f$-polynomial\/}
$$
f(x)=f_\De(x)=f_0+f_1x+\cdots+f_r x^r
$$
where $x$ is a variable.  
In the future we will continue the
practice of appending $\De$ in parentheses or as a subscript when we
wish to specify the complex, even if we do not do so in the
corresponding definition.

Another important invariant of $\De$ is its $h$-vector.  Define a
polynomial
$$
h(x):=(1-x)^r f\left(\frac{x}{1-x}\right)
=f_0(1-x)^r+f_1x(1-x)^r+\cdots+f_r x^r.
$$
Let $h_i$ be the coefficient of $x^i$ in $h(x)$ so that 
$h(x)=\sum_i h_i x^i$.  Then the {\it $h$-vector of $\De$\/} is
$$
\bh=(h_0,h_1,\ldots,h_r).
$$ 
It will sometimes be convenient to extend the range of
definition of the $f_i$ and $h_i$ by letting $f_i=h_i=0$ if $i<0$ or
$i>r$. 

Now suppose  that $G$ is a finite graph with vertices $V=V(G)$ and edges
$E=E(G)$.  We permit loops and multiple edges and will use the notation
$p=|V|$ and $q=|E|$.  We will also write $v\in G$ for $v\in V(G)$ and 
$e\in G$ for $e\in E(G)$ if it is clear from context whether we are
talking about the vertices or edges of $G$.
A {\it coloring\/} of $G$ is a function
$c:V\ra\{1,2,\ldots,\la\}$ and $c$ is {\it proper\/} if 
$c(u)\neq v(v)$ for all edges $uv\in E$.  Consider $G$'s {\it chromatic
polynomial}, $P(G)=P(G;\la)$, which is the number of such proper
colorings.  Note that if $G$ has a loop then $P(G;\la)=0$.
It is well known that if $G$ is a loopless then $P(G;\la)$ is a monic
polynomial 
of degree $p$ in $\la$ whose coefficients alternate in sign.
Writing
\beq
\label{Pf}
P(G;\la)=f_0\la^p-f_1\la^{p-1}+\cdots+(-1)^p f_p
\eeq
one can give the following interpretation to the coefficients $f_i$.

Let $\cC=\cC(G)$ denote the set of {\it cycles\/} of $G$ which will also be called
the set of {\it circuits\/}.  Suppose $G$ is {\it ordered\/} in that the edge set
$E$ has been given a linear ordering $e_1<e_2<\ldots<e_q$.  
Then each $C\in\cC$ gives 
rise to a {\it broken circuit\/} 
$$
\Cb=C-\min C
$$ 
where $\min C$ is the
smallest edge of $C$ in the linear ordering.
The {\it broken complex of $G$}, $\De(G)$, is the
family of all subsets of $E$ which do not contain a broken
circuit.  
It is easy to see that $\De(G)$ is a pure abstract simplicial complex.
Wilf~\cite{wil:wpc} was the first to consider this family of
sets as a complex.
In fact, $\De(G)$ is intimately connected with the chromatic
polynomial as can be seen in the following result which dates back to
Whitney~\cite{whi:lem}, although he did not state it in this form.
\bth[\cite{whi:lem}]
Let $P(G;\la)$ have coefficients $f_i$ as defined by~\ree{Pf}.  Then
$$
f_i=f_i(\De(G)),\ 0\le i\le p.\qqed
$$
\eth

One can think of the expansion~\ree{Pf} as being generated by a
sequence of deletions and contractions expressing $P(G;\la)$ as a
linear combination of chromatic polynomials of graphs with no edges.  One
could use chromatic polynomials of trees instead, or equivalently expand 
$P(G;\la)$ in terms of the basis $\{1\}\cup\{\la(\la-1)^i\ :\ i\ge0\}$
for the ring of polynomials in $\la$.  So define coefficients $h_i$ by
\beq 
\label{Ph}
P(G;\la)=h_0\la(\la-1)^{p-1}-h_1\la(\la-1)^{p-2}+\cdots+(-1)^p h_p.
\eeq
The next result follows easily from the previous theorem and the
definitions.
\bco
Define coefficients $h_i$ by~\ree{Ph}.  Then
$$
h_i=h_i(\De(G)),\ 0\le i\le p.\qqed
$$
\eco

Our goal is to give an explicit combinatorial description of the $h_i$
directly in terms of the broken circuits of the graph.  To do this, we
will need some 
machinery from the theory of Cohen-Macaulay rings.

\section{Cohen-Macaulay rings and monomial ideals}

Consider the polynomial ring $F[\bx]=F[x_1,x_2,\ldots,x_q]$ where
$F$ is a field and $\bx=\{x_1,x_2,\ldots,x_q\}$ is a set of variables.
If $E=\{e_1,e_2,\ldots,e_q\}$ then any $S\sbe E$ has 
corresponding monomial
$$
\bx^S=\prod_{e_i\in S} x_i.
$$
Now given any simplicial complex $\De$ on $E$ we form its  {\it
Stanley-Reisner ring}, $F(\De)$, by modding out by the non-faces of
$\De$, i.e., 
$$
F(\De)=F[\bx]/\spn{\bx^S\ :\ S\not\in\De}
$$
where $\spn{\cdot}$ denotes the ideal generated by the polynomials in the
brackets.  Note that since we are generating an ideal, it suffices
to consider the $\bx^S$ where $S$ is a minimal non-face of $\De$.

If $G$ is an ordered graph, then define
$$
F(G):=F(\De(G))=F[\bx]/\spn{\bx^\Cb\ :\ C\in\cC(G)}
$$
where we identify a (broken) circuit with its edge set.
This ring has a
{\it homogeneous system of parameters (hsop) of degree one\/},
i.e., a set of polynomials $\th_1,\ldots,\th_r\in F[\bx]$  which are homogeneous
of degree one and satisfy
\ben
\item $\th_1,\ldots,\th_r$ are algebraically independent, and
\item $F(G)/\spn{\th_1,\ldots,\th_r}$ is a finite dimensional vector space
  over $F$.
\een
Brown~\cite{bro:cpo}
gave an explicit construction of an hsop as follows.  To simplify
things, we will assume for now that $F=\bbZ_2$, the integers modulo
two.  In the last section, we will describe how to modify these ideas
so that they will work over an arbitrary field.

First note that if $G$ has blocks (maximal subgraphs having no
cutpoints) $G_1,G_2,\ldots,G_b$, then directly from the definitions we
have the ring isomorphism 
\beq
\label{iso}
F(G)\iso \otimes_i F(G_i).  
\eeq
So there is no
loss of generality in assuming that $G$ is a block and, in particular,
that $G$ is connected.  
Let $T$ be a spanning tree of $G$.  For each
edge $e\in T$, let $T'_e$ and $T''_e$ be the components of $T-e$.  So
$e$ defines a {\it fundamental cocircuit\/}
$$
D_e=D_e(G)=\{uv\in E(G)\ :\ u\in T'_e,\ v\in T''_e\}
$$ 
as well as a homogeneous degree one polynomial
\beq
\label{th}
\th_e=\sum_{e_i\in D_e} x_i.
\eeq

\thicklines
\setlength{\unitlength}{2pt}
\bfi
\bpi(80,70)(-10,-10)
\put(-10,30){\makebox(0,0){$G=$}}
\Gda
\Gad \Gdd \Ggd
\Gdg
\put(12,12){\makebox(0,0){$3$}}
\Gdaad
\put(48,12){\makebox(0,0){$7$}}
\Gdagd
\put(15,25){\makebox(0,0){$5$}}
\Gaddd
\put(12,48){\makebox(0,0){$1$}}
\Gaddg
\put(45,25){\makebox(0,0){$6$}}
\Gddgd
\put(33,43){\makebox(0,0){$2$}}
\Gdddg
\put(48,48){\makebox(0,0){$4$}}
\Ggddg
\epi
\hspace{50pt}
\bpi(80,70)(-10,-10)
\put(-10,30){\makebox(0,0){$T=$}}
\Gda
\Gad \Gdd \Ggd
\Gdg
\put(48,12){\makebox(0,0){$7$}}
\Gdagd
\put(15,25){\makebox(0,0){$5$}}
\Gaddd
\put(45,25){\makebox(0,0){$6$}}
\Gddgd
\put(48,48){\makebox(0,0){$4$}}
\Ggddg
\epi
\capt{A graph $G$ and spanning tree $T$}\label{GT}
\efi

Since this construction will be crucial, we illustrate it with an
example.  Consider the graph $G$ and its spanning tree $T$ given in
Figure~\ref{GT}.  For simplicity we have labeled the edges
$1,2,\ldots,7$ rather than $e_1,e_2,\ldots,e_7$.  Then we have
\bea
\th_4&=&x_4+x_1+x_2,\\
\th_5&=&x_5+x_1+x_3,\\
\th_6&=&x_6+x_1+x_2+x_3,\\
\th_7&=&x_7+x_3.\\
\eea
For any graph $G$ we have the following result.
\bth[\cite{bro:cpo}]
If $G$ is a connected graph and $T$ a spanning tree then 
the set of polynomials defined by~\ree{th} for $e\in T$ is an hsop for
$\bbZ_2(G)$.\qqed 
\eth

Continuing with the general development, let  $\Mon(\bx)=\Mon(q)$
denote the set 
of monomials in $F[\bx]=F[x_1,x_2,\ldots,x_q]$.  When it will do no
harm, we will not distinguish between these monomials considered as
elements of $F[\bx]$ or considered as elements of some quotient of the
polynomial ring.  A subset $L\subseteq\Mon(q)$ is a {\it lower order
ideal\/} (or {\it down set\/})
if whenever $m\in L$ and $n\in\Mon(q)$ divides $m$, then $n\in L$.
Similarly, $U \subseteq\Mon(q)$ is an {\it upper order
ideal\/} (or {\it filter\/}) if whenever $m\in L$ and $n\in\Mon(q)$ is
divisible by $m$, then $n\in L$.  Note that $U$ is an upper order ideal
if and only if $\Mon(q)-U$ is a lower order ideal.
If $S\subseteq\Mon(q)$ then the {\it lower and upper order ideals
generated by $S$\/} are
\bea
L(S)&=&\{n\in\Mon(q)\ :\ \mbox{$n$ divides $m$ for some $m\in S$}\},\\
U(S)&=&\{n\in\Mon(q)\ :\ \mbox{$n$ is divisible by $m$ for some $m\in S$}\}.
\eea
Macaulay~\cite{mac:pet} showed  that after modding out by an hsop, one can always
find a basis of monomials which forms a lower order ideal.  And
Stanley~\cite{sta:ubc} connected such a basis with the $h$-vector.  
\bth[\cite{mac:pet,sta:ubc}]
\label{mac}
Suppose that $I$ is an ideal of $F[\bx]$ and that
$\th_1,\ldots,\th_r$ be an hsop for $F[\bx]/I$.  Then
the ring
$$
R=\frac{F[\bx]}{I+\spn{\th_1,\ldots,\th_r}}
$$
has a basis $L$ which is a lower order ideal of monomials.  

Suppose further that $F[\bx]/I$ is Cohen-Macaulay and $F[\bx]/I\iso
F(\De)$ for some simplicial complex 
$\De$ with $h$-vector $\bh=(h_0,\ldots,h_r)$.  Then
$$
h_i = \mbox{number of monomials of total degree $i$ in $L$}.\qqed
$$
\eth

Now consider a graph $G$ with a spanning tree $T$ and define $I(G)$ to
be the ideal of $F[\bx]$ generated by the monomials $\bx^\Cb$ for
$C\in\cC(G)$.  We wish to give an
explicit basis for the ring
$$
R(G)=\frac{F[\bx]}{I(G)+\spn{\th_e\ :\ e\in T}}
$$ 
which is a lower order ideal of monomials.  First, however,  we 
wish to show that we have a basis inside $\Mon(\by)$ for a subset
$\by$ of $\bx$.

An ordering $e_1<e_2<\ldots<e_q$ of
$E(G)$ will be called {\it standard\/} if the last $p-1$ edges in the
order form a tree.  From now on we will assume that all our orderings
are standard and take our spanning tree $T=T(G)$ to be the one determined
the last edges in the order.  It will also be convenient to denote the
number of edges not in $T$ by $k=q-p+1$.  We will show that we that
our basis can be taken in $\Mon(\by)$ where $\by=\{x_1,x_2,\ldots,x_k\}$.

We now return to working over $\bbZ_2$.
Suppose $k<j\le q$ and write $D_j$ for
$D_{e_j}$ and $\th_j$ for $\th_{e_j}$.  Then since $\th_j=0$ in $R(G)$
we have
\beq
\label{xj}
x_j=\sum_{e_i\in (D_j-e_j)} x_i.
\eeq
where $x_i\in\by$ for all $x_i$ appearing in the sum.  For each $C\in\cC$
let $p_\Cb=p_\Cb(\by)$ be the polynomial obtained from $x^\Cb$ by
substituting in 
the sum in equation~\ree{xj} for $x_j$ for each $j>k$.
Consider the ideal
$$
J=J(G)=\spn{p_\Cb\ :\ C\in\cC}
$$
We immediately have the following result.
\bpr
If $G$ is a connected graph and $F=\bbZ_2$ then
$$
R(G)\iso \frac{\bbZ_2[\by]}{J(G)}.\qqed
$$
\epr

Returning to our running example, we convert the list of circuits in
$G$ into polynomials using the equations for
$\th_4,\ldots,\th_7$.
$$
\barr{lll}
C_1=\{1,4,5,6\},\qquad  &\bx^{\ol{C_1}}=x_4x_5x_6,\qquad   
      &p_{\ol{C_1}}=(x_1+x_2)(x_1+x_3)(x_1+x_2+x_3),\\
C_2=\{2,4,6\}           &\bx^{\ol{C_2}}=x_4x_6,         
      &p_{\ol{C_2}}=(x_1+x_2)(x_1+x_2+x_3),\\
C_3=\{3,5,6,7\}         &\bx^{\ol{C_3}}=x_5x_6x_7,      
      &p_{\ol{C_3}}= x_3 (x_1+x_3)(x_1+x_2+x_3),\\
C_4=\{1,2,5\}           &\bx^{\ol{C_4}}=x_2x_5,         
      &p_{\ol{C_4}}= x_2 (x_1+x_3),\\
C_5=\{1,3,4,7\}         &\bx^{\ol{C_5}}=x_3x_4x_7,      
      &p_{\ol{C_5}}= x_3^2 (x_1+x_2),\\
C_6=\{2,3,4,5,7\}       &\bx^{\ol{C_6}}=x_3x_4x_5x_7,      
      &p_{\ol{C_6}}= x_3^2(x_1+x_2)(x_1+x_3)\\
C_7=\{1,2,3,6,7\}       &\bx^{\ol{C_7}}=x_2x_3x_6x_7,      
      &p_{\ol{C_7}}= x_2 x_3^2(x_1+x_2+x_3).\\
\earr
$$

We will now pick a specific monomial  $m_\Cb$ from 
each $p_\Cb$ and these will be used to define the lower order ideal of
monomials being sought.
For $1\le i\le k$, the graph $T+e_i$ has a
unique circuit $C_i$ and these circuits will be called {\it fundamental}.
We label the nonfundamental circuits in some order as $C_i$ for $i>k$.
Also define
$$
d_i=\case{i}{if $i\le k$,}{\min \{j\ :\ e_j\in D_i\}}{if $i>k$.}
$$
Now let
$$
m_{\Cb_i}=\case{x_i^{|\Cb_i|}}{if $i\le k$,}
{\rule{0pt}{30pt}\dil\prod_{e_j\in\Cb_i} x_{d_j}}{if $i>k$.}
$$
It is easy to see from the definitions that $m_\Cb$ is indeed a term
in the polynomial $p_\Cb$.
Finally, define upper and lower order ideals
$$
U(G)=U(m_\Cb\ :\ C\in\cC(G))\qmq{and} L(G)=\Mon(k)-U(G).
$$
Note that all these quantities depend on the ordering imposed on the
edges and not just on the graph itself, even though our notation does
not reflect that.  It is $L(G)$ which will be our candidate as a
monomial basis for $R(G)$

Continuing with our example, $C_1$, $C_2$, and $C_3$ are fundamental
with 4, 3, and 4 edges (respectively) and so
$$
m_{\Cb_1}=x_1^3,\quad m_{\Cb_2}=x_2^2,\quad m_{\Cb_3}=x_3^3.
$$
The monomials $m_{\Cb}$ for the other four circuits are obtained by taking the
variable of smallest subscript in each factor of the corresponding
$p_{\Cb}$, so
$$
m_{\Cb_4}=x_1 x_2,\quad m_{\Cb_5}=x_1 x_3^2,\quad 
m_{\Cb_6}=x_1^2 x_3^2,\quad m_{\Cb_6}=x_1 x_2 x_3^2.
$$
Thus $R(G)$ should have as basis
\bea
L(G)
&=&
\Mon(3)-U(x_1^3,\ x_2^2,\ x_3^3,\ x_1 x_2,\  x_1 x_3^2,\  x_1^2x_3^3,\ x_1 x_2 x_3^2)\\
&=&
\{1,\ x_1,\ x_2,\ x_3,\ x_1^2,\ x_3^2,\ x_1x_3,\ x_2x_3,\ x_1^2 x_3,\ x_2 x_3^2 \}.
\eea
and this can be verified directly.

A graph for which there is an ordering of $E$ such that $L(G)$ is a
basis for $R(G)$ will be said to have a {\it
no broken circuit basis\/} or {\it NBC basis\/}.  
To outline the rest of the paper,
in the next section we will prove a
general theorem about when a graph has an NBC basis.
In Section~\ref{tf} we will apply these results
to show that two infinite
families of graphs do indeed have NBC bases.  
Section~\ref{ub} will be devoted to 
giving an upper bound for the number of acyclic orientations for a graph with an NBC basis.
We also compare this bound to others in the literature.
We end with some comments and open problems.  This will include a 
conjecture that every graph $G$ has ordering which produces an NBC
basis for $R(G)$, as well as a proposed line of attack on this idea.

\section{Graphs with NBC bases}
\label{gwb}

One way to show that a graph has an NBC basis would be to use
induction.  Since 
the chromatic polynomial is involved, this would entail deletion and
contraction.  If $e\in E(G)$ then let $G\sm e$ and $G/e$ denote $G$
with $e$ deleted and with $e$ contracted, respectively.
Since we are permitting loops and multiedges, both $G\sm e$ and $G/e$
will have exactly one less edge than $G$.
An elementary fact about the chromatic polynomial is that
$$
P(G;\la)=P(G\sm e;\la)-P(G/e;\la).
$$
Using this equation and~\ree{Ph} we easily obtain 
the following proposition.
\bpr
Let $G$ be a graph and $e\in E(G)$.  Then for all $i\ge0$ we have
$$
h_i(G)=h_i(G\sm e)+h_{i-1}(G/e).\qqed
$$
\epr

If we choose $e\in T$ then $T/e$ is a spanning tree of $G/e$ but 
$T\sm e$ is no longer a tree.  If, on the other hand, we choose 
$e\not\in T$ then $T$ is still a spanning tree of $G\sm e$ but $T$ is
no longer a tree in $G/e$.  However, we can get around these
difficulties if $G$ has a vertex $w$ with $\deg w =2$ where $\deg w$,
the {\it degree of $w$}, is the number of edges containing $w$.

As noted before, it does no harm to restrict our attention to graphs
$G$ which are blocks so that $G\sm e$ and $G/e$ are connected for all
$e\in E$.   We will say that a standard ordering $e_1<e_2<\ldots<e_q$ on $G$
imposes the {\it induced ordering\/} $e_1<e_2<e_3<\ldots<e_{q-1}$ on
$G\sm e_q$ and  on $G/e_q$.  Now suppose that $G$ has a
vertex $w$ with $\deg w=2$ and that $e_k,e_q$ are the two edges
containing $w$.  Then if the ordering on $G$ is standard, so too
will be the induced orderings on $G\sm e_q$ and $G/e_q$.   
Our primary tool for showing that certain graphs have NBC
bases will be the following theorem.  Note that an example which
illustrates the proof of this result follows the demonstration, so the
reader may wish to read both in parallel.

\bth
\label{deg2}
Let $G$ be a block with a standard ordering $e_1<e_2<\ldots<e_q$.
Suppose $G$ has a vertex $w$ of degree two such that the edges
containing $w$ are $e_k$ and $e_q$.  If $R(G\sm e_q)$ and $R(G/e_q)$
have NBC bases in their induced standard orderings, then so does
$R(G)$.
\eth
\pf\
Let $\uplus$ denote disjoint union and if $S\sbe\Mon(k)$ and
$m\in\Mon(k)$ then let $mS=\{mn\ :\ n\in S\}$.
We first show that
\beq
\label{L(G)}
L(G)=L(G\sm e_q)\uplus x_k L(G/e_q)
\eeq
so that by our assumptions about 
$R(G\sm e_q)$ and $R(G/e_q)$ and the previous proposition (summed
over all $i$) we have
$$
|L(G)|=|L(G\sm e_q)|+|L(G/e_q)|=\dim R(G\sm e_q)+\dim R(G/e_q)=\dim R(G)
$$
where dimension is being taken over the field $\bbZ_2$.

Consider $G\sm e_q$.  Note that $e_k$ is in the tree for $G\sm e_q$
and so the basis for $R(G\sm e_q)$ will be in $\Mon(k-1)$.
Also, from our assumptions on $w$, $e_k$ is the only
edge of $G\sm e_q$ containing $w$.  So $C$ is a circuit of $G\sm e_q$ if
and only if $C$ is a circuit of $G$ not containing $e_k$.  It follows
that $x_k$ is never a factor of $\bx^\Cb$ for such $C$.  It also
follows that for $e_j\in T(G\sm e_q)$, $e_j\neq e_k$, we have
$D_j(G\sm e_q)=D_j(G)-e_k$. And both of these sets have the same
minimum since $e_k$ is the edge of largest index outside the tree for
$G$.  Thus the generators for $J(G\sm e_q)$ are 
obtained from those for $J(G)$ by setting $x_k=0$ wherever it
appears.  So the monomials in $U(G\sm e_q)$ are precisely those in
$U(G)$ which do not have $x_k$ as a factor.  Hence $L(G\sm e_q)$
consists of the monomials in $L(G)$ which do not have $x_k$ as a
factor.

Now consider $G/e_q$.  The circuits of $G$ are in bijection with the
circuits of $G/e_q$:  If $C\in\cC(G)$ contains $e_q$ then it 
corresponds to the circuit $C/e_q$ of $G/e_q$, while if $C$ does not
contain $e_q$ then it is
also a circuit of $G/e_q$ itself.  We will call the former
circuits (in both $G$ and $G/e_q$) {\it type I\/}, and the latter {\it
type II}.  Note that because of the assumptions on $w$, the type I and
type II circuits can also be characterized as those which do and do not
contain $e_k$, respectively.  Since $e_q$ is the only edge of $T(G)$
containing $w$, we have $D_j(G/e_q)=D_j(G)$ for each 
$e_j\in T(G/e_q)$.  Thus, using $\pt_\Cb$ to denote the generators of
$J(G/e_q)$, 
$$
p_\Cb=\case{x_k \pt_{\ol{C/e_q}}}{if $C$ is of type I,}
{\pt_\Cb}{if $C$ is of type II,}
$$
where the polynomials for the type II circuits have no factor of $x_k$.
Since $e_k$ has the largest index outside $T(G)$, the same relation
holds between the corresponding generators of $U(G)$ and of
$U(G/e_q)$, i.e., $m_\Cb= x_k {\tilde{m}}_{\ol{C/e_q}}$ or 
${\tilde{m}}_\Cb$ depending on whether $\Cb$ is type I or type II
(respectively), where the tilde indicates the the quantity is being
calculated in $C/e_q$.

Now  one sees that
$x_k L(G/e_q)$ consists precisely of the
monomials in $L(G)$ which have a factor of $x_k$:
Suppose that we have a monomial of $L(G)$ divisible by $x_k$.  Then it can be
written as $x_k m$ for some $m\in\Mon(k)$.  Since $x_k m$ is not
divisible by any type I generator of $U(G)$, and all such generators
have the form $x{\tilde{m}}_1$ for some type I generator ${\tilde{m}}_1$ of $U(G/e_q)$,
we see that $m$ is not divisible by ${\tilde{m}}_1$ for all type I generators
of $U(G/e_q)$.  Also, $x_k m$ is not divisible by any type II
generator ${\tilde{m}}_2$ of $U(G)$, and all such generators do not have $x_k$
as a factor, so $m$ is not divisible by any type II generator ${\tilde{m}}_2$ of
$U(G/e_q)$.  So $m\in L(G/e_q)$ and $x_k m \in x_k L(G/e_q)$.  The
proof of the converse inclusion is similar.

Since $L(G)$ is
clearly the disjoint union of its monomials with a factor of $x_k$ and
its monomials without a factor of $x_k$, we are done with the
demonstration of~\ree{L(G)}.
So we have proven that $L(G)$ contains $\dim R(G)$
monomials, and thus it will
suffice to show that these monomials span $R(G)$.  For that, it
suffices to show that $L(G)$ spans $U(G)$.  So take $m\in U(G)$.
Suppose first that $x_k$ is a factor of $m$ so that $m=x_kn$ for some
monomial $n$.  Then from our work in the previous paragraph we see
that $n\in U(G/e_q)$.  So by our assumption about $R(G/e_q)$, we can
write
\beq
\label{n}
n=\sum_{l\in L(G/e_q)} a_l l+ p
\eeq
where the $a_l$ are constants and $p\in J(G/e_q)$.
But $x_k l\in L(G)$ for $l\in L(G/e_q)$, and $x_k p\in J(G)$ for 
$p\in J(G/e_q)$ since this is true for each of the generators of 
$J(G/e_q)$.  So multiplying~\ree{n} by $x_k$ expresses
$m=x_kn$ as a linear combination of elements of $L(G)$ modulo $J(G)$
as desired.

Now suppose that $x_k$ is not a factor of $m$.  Then by our previous
results concerning $G\sm e_q$ we see that $m\in U(G\sm e_q)$.  So
by our assumption about $R(G\sm e_q)$, we can write
\beq
\label{m}
m=\sum_{l\in L(G\setminus e_q)} a_l l + p
\eeq
where the $a_l$ are constants and $p\in J(G\sm e_q)$.
Now, as shown above, $l\in L(G\sm e_q)$ implies $l\in L(G)$.
Furthermore, there must be a $p'\in J(G)$ such that
$p'(x_1,\ldots,x_{k-1},0)=p$.  It follows that $p'=p+x_k p''$ for some
$p''\in F[\by]$.  But, from the previous paragraph, we have that
$x_kp''$ is spanned by $L(G)$ modulo $J(G)$.  So substituting
$p=p'-x_kp''$ into~\ree{m} we have expressed $m$ as a linear
combination of elements of $L(G)$ modulo $J(G)$.  Hence every monomial
is in the span of $L(G)$ and we are done.\Qqed

\thicklines
\setlength{\unitlength}{2pt}
\bfi
\bpi(100,70)(-30,0)
\put(-20,30){\makebox(0,0){$G\sm e_7=$}}
\Gda
\Gad \Gdd \Ggd
\Gdg
\put(12,12){\makebox(0,0){$3$}}
\Gdaad
\put(15,25){\makebox(0,0){$5$}}
\Gaddd
\put(12,48){\makebox(0,0){$1$}}
\Gaddg
\put(45,25){\makebox(0,0){$6$}}
\Gddgd
\put(33,43){\makebox(0,0){$2$}}
\Gdddg
\put(48,48){\makebox(0,0){$4$}}
\Ggddg
\epi
\hspace{50pt}
\bpi(100,70)(-30,0)
\put(-20,30){\makebox(0,0){$T(G\sm e_7)=$}}
\Gda
\Gad \Gdd \Ggd
\Gdg
\put(12,12){\makebox(0,0){$3$}}
\Gdaad
\put(15,25){\makebox(0,0){$5$}}
\Gaddd
\put(45,25){\makebox(0,0){$6$}}
\Gddgd
\put(48,48){\makebox(0,0){$4$}}
\Ggddg
\epi
\capt{The graph $G\sm e_7$ and spanning tree $T(G\sm e_7)$}\label{Gd}
\efi

Returning to our example graph (which satisfies the conditions of the
previous theorem), $G\sm e_7$ and the tree for the induced order are
shown in Figure~\ref{Gd}.  The relevant sets are  
\bea
\{\th_e\}
&=&
\{x_3,\ x_4+x_1+x_2,\ x_5+x_1,\ x_6+x_1+x_2\},
\\
\{C\}
&=&
\{\{1,4,5,6\},\ \{2,4,6\},\ \{1,2,5\}\},
\\
\{\bx^\Cb\}
&=&
\{x_4 x_5 x_6,\ x_4 x_6,\ x_2 x_5\},
\\
\{p_\Cb\}
&=&
\{x_1(x_1+x_2)^2,\ (x_1+x_2)^2,\ x_1 x_2\},
\\
U(G\sm e_7)
&=&
U(x_1^3,\ x_2^2,\ x_1 x_2),
\\
L(G\sm e_7)
&=&
\Mon(2)-U(G\sm e_7)
\\
&=&
\{1,\ x_1,\ x_2,\ x_1^2\}.
\eea

Making the same computations in $G/e_7$ yields
\bea
\{\th_e\}
&\hs{-10pt} = \hs{-10pt}&
\{x_4+x_1+x_2,\ x_5+x_1+x_3,\ x_6+x_1+x_2+x_3\},
\\
\{C\}
&\hs{-10pt} = \hs{-10pt}&
\{\{1,4,5,6\},\ \{2,4,6\},\ \{3,5,6\},\ \{1,2,5\},\ \{1,3,4\},\
\{2,3,4,5\},\ \{1,2,3,6\}\},
\\
\{\bx^\Cb\}
&\hs{-10pt} = \hs{-10pt}&
\{x_4 x_5 x_6,\ x_4 x_6,\ x_5 x_6,\ x_2 x_5,\ x_3 x_4,\ x_3x_4x_5,\ x_2x_3x_6\},
\\
\{p_\Cb\}
&\hs{-10pt} = \hs{-10pt}&
\{(x_1+x_2)(x_1+x_3)(x_1+x_2+x_3),\ 
(x_1+x_2)(x_1+x_2+x_3), (x_1+x_3)(x_1+x_2+x_3),\
\\
&& 
\hs{5pt} x_2(x_1+x_3),\ x_3(x_1+x_2),\ x_3(x_1+x_2)(x_1+x_3),\ x_2 x_3(x_1+x_2+x_3)\},
\\
U(G/e_7)
&\hs{-10pt} = \hs{-10pt}&
U(x_1^3,\ x_2^2,\ x_3^2,\ x_1 x_2,\ x_1 x_3,\ x_1^2 x_3,\ x_1 x_2 x_3),
\\
L(G/e_7)
&\hs{-10pt} = \hs{-10pt}&
\Mon(3)-U(G/e_7)
\\
&\hs{-10pt} = \hs{-10pt}&
\{1,\ x_1,\ x_2,\ x_3,\ x_1^2,\ x_2 x_3 \}.
\eea
Note that we have $L(G)=L(G\sm e_7)\uplus x_3 L(G/e_7)$.

\thicklines
\setlength{\unitlength}{2pt}
\bfi
\bpi(100,70)(-30,0)
\put(-20,30){\makebox(0,0){$G/e_7=$}}
\Gad \Gdd \Ggd
\Gdg
\put(30,5){\makebox(0,0){$3$}}
\put(30,30){\oval(60,40)[b]}
\put(15,25){\makebox(0,0){$5$}}
\Gaddd
\put(12,48){\makebox(0,0){$1$}}
\Gaddg
\put(45,25){\makebox(0,0){$6$}}
\Gddgd
\put(33,43){\makebox(0,0){$2$}}
\Gdddg
\put(48,48){\makebox(0,0){$4$}}
\Ggddg
\epi
\hspace{50pt}
\bpi(100,70)(-30,0)
\put(-20,30){\makebox(0,0){$T(G/e_7)=$}}
\Gad \Gdd \Ggd
\Gdg
\put(15,25){\makebox(0,0){$5$}}
\Gaddd
\put(45,25){\makebox(0,0){$6$}}
\Gddgd
\put(48,48){\makebox(0,0){$4$}}
\Ggddg
\epi
\capt{The graph $G/e_7$ and spanning tree $T(G/e_7)$}\label{Gc}
\efi

\section{Two families}
\label{tf}

We will now consider two families of graphs and prove that they have
NBC bases.  They are called (generalized) theta and phi
graphs.

A {\it (generalized) theta graph\/} consists of two vertices $u,v$
together with $t$ internally-disjoint 
$u\hor v$ paths $P',P'',\ldots,P^{(t)}$.  Note that we are not insisting
that $t=3$ as is usually done for theta graphs.   
To show that such a
graph has an NBC basis, we need to label its edges so that
$e_1<e_2<\ldots<e_q$ is a standard order.
First label all the edges in paths of length one with
$e_k,e_{k-1},\ldots,e_{l+1}$ for some $l\le k$.  Now take any remaining
path of length at least two and label its edges, starting from the one
containing $u$, as $e_l,e_q,e_{q-1},e_{q-2},\ldots,e_{r+1}$ for some
$r$.  Now take another path of length at least two (if any) and label
its edges $e_{l-1},e_r,e_{r-1},e_{r-2},\ldots,e_{s+1}$.  Continue in
this way until all 
the edges have been labeled.  Note that this labeling does produce
a standard ordering and will be called a {\it theta labeling}.

\bth  If $G$ is a (generalized) theta graph with a theta labeling then
$G$ has an NBC basis.
\eth
\pf\
We will induct on the number of edges of $G$.  If $G$ is a single path
or if all paths are of length one, then the result is easy to verify.
So we may assume that $G$ is a block and has at least one $u\hor v$
path $P'$ of length two or greater

Let $w$ be the vertex on $P'$ adjacent to $u$.  Then we have set
things up so that $w$ satisfies the hypotheses of Theorem~\ref{deg2}
except that $w$ is adjacent to $e_q$ and $e_l$ for some $l\le k$, not
necessarily $e_k$ itself.  But the reason we chose $e_k$ in the
proof of the theorem was because $k$ was the largest index outside
$T(G)$.  This guaranteed that for each circuit $C$, the monomials
$m_\Cb$ picked from the $p_\Cb$ in $G$, $G\sm e_q$, and $G/e_q$ would be
related in the correct way.  And the reason for this was that given
any edge $e$ of $G$ which was both in a circuit and in $T(G)$, the
cocircuit $D_e$ would contain an edge of index smaller than $k$ and so 
$x_k$ would not be picked from that factor.  But because of the way we
have chosen to label the $u\hor v$ paths of length one, the preceeding
statements also hold if one replaces $e_k$ by $e_l$ everywhere.  So
this change in index does no harm and will permit us to use induction,
as a theta labeling of $G$ will induce theta labelings of $G\sm e_q$
and $G/e_q$.

Now consider $G\setm e_q$.  
This is not a theta graph in general.  
But the induced labeling on $G\sm e_q$ is a theta
labelling if we ignore the other edges on $P'$.  This does not cause
any problems since each of these edges is now a block and so does not
contribute anything to $F(G)$ by~\ree{iso} and the fact that 
$R(e)\iso F$ for any edge $e$.  Hence, by induction, $R(G\sm e_q)$ has an
NBC basis.

Now look at $G/e_q$.  This is still a theta graph and, since $P'$ has
length at least two, its induced
labeling is a theta labeling.  So, by induction, $R(G/e_q)$ has an NBC
basis.  Hence all the hypotheses of Theorem~\ref{deg2} are
satisfied and $G$ has an NBC basis, completing our proof.
\Qqed

As a special case of the previous result, we obtain the following.
\bco
The complete bipartite graph $K_{2,t}$ with a theta labeling has an
NBC basis.\qqed
\eco

Rather than thinking of theta graphs as unions of paths, one could
consider them as a set of cycles joined in parallel.  We will now
define a family of graphs which can be thought of as joining cycles in
series.  Suppose we are given $t$ cycles $C',C'',\ldots,C^{(t)}$ all
of length at least two, and in
each $C^{(i)}$ we are given a pair of distinguished edges
$e^{(i)},f^{(i)}$.  Then 
the associated {\it phi graph\/} is obtained by identifying $f^{(i)}$ with
$e^{(i+1)}$ for $1\le i< t$.  
For example, if we let $P_p$ denote the
path on $p$ vertices then the cross product $P_2\times P_t$ is a phi
graph where all the cycles have length four.  
(It is because of the shape of $P_2\times P_3$ that we call these phi graphs.)

Again, we will need a specific labeling for our phi graphs.  Label
edge $e^{(i)}$ with $e_{k-i+1}$, $1\le i\le t$.  Now label the remaining
edges of $C^{(1)}$ as follows.  We have
$C^{(1)}-e^{(1)}-f^{(1)}=P\uplus Q$ where $P, Q$ are paths.  Label the
edges along $P$ (if any) starting with the one adjacent to $e^{(1)}$
with $e_q,e_{q-1},\ldots,e_{r+1}$.  Now do the same along $Q$ using
the labels $e_r,e_{r-1},\ldots,e_s$.  Continue in like manner to label
the rest of the cycles.  (When one gets to the last one, there will be
only one path to label.)  Call this a {\it phi labeling\/} of the
graph.

Before proving that a phi graph has an NBC basis, we will
need a lemma to take care of the special case when the first cycle has
length two, so that attaching it to the second cycle creates an edge of
multiplicity two.  
Let $G$ be a connected graph with standard ordering
$e_1<e_2<e_3<\ldots<e_q$ where $e_k$ and $e_{k-1}$ have the same
endpoints.  Let $G\sm e_k$ have the induced ordering
$e_1<\ldots<\eh_k<\ldots<e_q$ where the hat indicates that $e_k$ has
been removed.  Note that the induced ordering is standard.
Then the corresponding rings are related in the manner in which one
would expect given that the chromatic polynomials do not change.
\ble
Suppose that $G$ has a standard ordering such that $e_k$ and $e_{k-1}$
have the same endpoints.  If $G\sm e_k$ is given the induced ordering
above then $R(G)\iso R(G\sm e_k)$.
\ele
\pf
Directly from the definitions one sees that one obtains the generators
for $J(G)$ from those for
$J(G\sm e_k)$ by substituting $x_{k-1}+x_k$ everywhere one has an $x_{k-1}$.  The
additional cycle made by $e_{k-1},e_k$ also sets $x_k=0$ in the quotient $R(G)$.
Hence the isomorphism.
\qqed

\bth
If $G$ is a phi graph with a phi labeling then
$G$ has an NBC basis.
\eth
\pf
Again, we induct on the number of edges in $G$.  The case of a single cycle
is easy to do (and appears in~\cite{bro:cpo}).  So suppose we have at
least two cycles. If $C'$ has length two, then its phi labeling is
exactly the type considered in the previous lemma.  So 
$R(G)\iso R(G\sm e_k)$ where the latter graph has a phi labeling and
fewer edges.  So we are done in this case.

If $C'$ has length at least three, then a deletion-contraction
argument similar to the one used for theta graphs will provide a
proof.  We leave the details of the demonstration to the reader.
\qqed

\bco
The graph $P_2\times P_t$ with a phi labeling has an NBC basis.\Qqed
\eco

\section{Upper Bounds}
\label{ub}

If graph $G=(V,E)$ has an NBC basis, then we can use this fact to give a
simple upper bound on its $h$-vector.  (Lower bounds for $h$-vectors
of various types of complexes have been given
by Swartz~\cite{swa:lbh}.) This, in turn, bounds the values
of the chromatic polynomial $P(G;\la)$ at negative integers since then
all terms in the expansion~\ree{Ph} have the same sign.  In
particular, this gives an upper bound on $\al(G)$, the number of
acyclic orientations of $G$, because of a famous theorem of
Stanley~\cite{sta:aog} which states that
$$
\al(G)=(-1)^p P(G;-1)
$$
where, as usual, $p=|V|$.  To see why one could only expect to bound
these quantities, rather than obtaining their exact values, we need to
say a few words about the theory of \#P problems which was introduced by
Valiant~\cite{val:ccp,val:cer}.  

If $A$ and $B$ are two problems then we say that $A$ is {\it
polynomially reducible\/} to $B$ if it is possible, given a subroutine
to solve $B$, to solve $A$ in polynomial time, where we count calls to
the subroutine for $B$ as a single step.
The {\it class \#P\/} consists of those enumeration problems where the
structures being counted can be recognized in polynomial time.  In
other words, there is an algorithm which is polynomial in the size
of the input problem that can verify whether a given structure should
be included in the count.  So the class \#P is to enumeration
problems as the class NP is to decision problems.  An enumeration
problem is {\it \#P-complete\/} if any problem in \#P is polynomially
reducible to it.  So the \#P-complete problems are the hardest in
\#P.  

Linial~\cite{lin:hep} first showed that computing $\al(G)$ is
\#P-complete.  Jaeger, Vertigan, and Welsh~\cite{jvw:ccj} derived more
general results about computing the Tutte polynomial of a matroid
which imply that computing $P(G;\la)$ is \#P-complete for all but
nine special values of $\la$.  

The case $\la=-1$ has attracted special interest because $\log\al(G)$
is a lower bound on the computational complexity of certain decision
and sorting problems, see for example the paper of Goddard, Kenyon,
King, and Schulman~\cite{gkks:ora}.  Obviously the number of acyclic
orientations of $G$ is bounded above by the total number of
orientations, giving
$$
\al(G)\le 2^q
$$
where $q=|E|$.  Fredman (whose work is reported in a paper of Graham,
Yao, and Yao~\cite[Section 7]{gyy:ibw}), and independently Manber and
Tompa~\cite{mt:enh} gave the first nonobvious upper bound for $\al(G)$ as
$$
\al(G)\le \prod_{v\in V} (\deg v+1),
$$ 
where, as usual, $\deg v$ is the degree of vertex $v$.  This bound was improved by
Kahale and Schulman~\cite{ks:bcp} as follows.

Given a graph $G$, consider its {\it cone\/}, $G^*$, obtained by
adding a new vertex adjacent to every vertex of $G$.  Then Kahale and
Schulman show that $\al(G)$ is at most the number of spanning trees of
$G^*$.  Using the Matrix-Tree Theorem, this bound can be expressed as
a determinant.  Since the determinant itself could be costly to
compute, they give an upper bound for its value.
\bth[\cite{ks:bcp}]
\label{KS}
We have the upper bound
\beq
\label{be}
\al(G)\le\prod_{v\in V} (\deg v+1)
\prod_{uw\in E} \exp\frac{-1}{2(\deg u+1)(\deg w+1)}
\stackrel{\rm def}{=}\be(G).\qqed
\eeq
\eth

Now suppose that $G$ has an NBC basis 
$\Mon(k)-U(m_\Cb\ :\ C\in\cC(G))$.  If we remove the upper order
ideal generated by just the fundamental circuits, then we will get a
spanning set for the quotient which can be used to bound the
$h$-vector from above.  Furthermore, each of these monomials  
has the simple form
$$
m_{\Cb_i}=x_i^{|C_i|-1}.
$$
So by Theorem~\ref{mac} and equation~\ree{Ph}, we have proved the
following result, where we use  
$L_d(S)$ to denote the set of monomials in the lower order ideal
$L(S)$ which have total degree $d$.
\bth
If $G$ has an NBC basis with fundamental circuits $C_1,\ldots,C_k$ then,
for $d\ge0$,
\beq
\label{ld}
h_d(G)\le \left|L_d\left(x_1^{|C_1|-2}\cdots x_k^{|C_k|-2}\right)\right|
\stackrel{\rm def}{=}l_d(G).
\eeq
Furthermore
\beq
\label{ga}
\al(G)\le\sum_{d=0}^{p-1} l_d(G) 2^{p-d-1}
\stackrel{\rm def}{=}\ga(G).\qqed
\eeq
\eth
We note that it is an easy exercise to show that
\beq
\label{ch}
l_d(G)\le|L_d(\Mon(k))|={d+k-1\choose k-1}.
\eeq
If one wishes, one can calculate the exact values of the $l_d(G)$
using the Principle of Inclusion-Exclusion (see Stanley's
text~\cite[Chapter 2]{sta:ec1}).

We will now compare the bounds $\be(G)$ and $\ga(G)$ for certain theta and phi
graphs.  When possible, we will compare the $\ga$ bound with the
actual number of acyclic orientations.  Of course, from a practical
viewpoint, it is unnecessary to use a bound when the exact value is
known.  But this will give some sense of how close $\ga$ is to the
truth.

We keep the conventions of the previous section.
Define $\Th_{n,t}$ to be the theta graph consisting of $t$
paths of length $n$ with their endpoints identified to form the
special vertices $u$ and $v$.  There is an interesting change in the
behaviour of the $\ga$ bound depending on whether $n$ is held fixed and
$t$ varies, or vice-versa.
\bth
As $n\ra\infty$ we have
$$
\ga(\Th_{n,3})\sim\al(\Th_{n,3}).
$$
As $t\ra\infty$ we have
$$
\be(\Th_{2,t})=o(\ga(\Th_{2,t})).
$$
\eth
\pf\
First consider $\Th_{n,3}$ where $p=3n-1$ and $q=3n$.  Since this
graph only has 3 circuits, it is easy to use Inclusion-Exclusion to
calculate $\al(G)$, from which one sees that the count is asymptotic
to the first term
$$
\al(G)\sim 2^q=2^{3n}.
$$

To compute $\ga$, first note that from~\ree{ld} and~\ree{ch} we have
$$
h_d(\Th_{n,3})\le l_d(x_1^{2n-2}\ x_2^{2n-2})\le d+1.
$$
Plugging this bound into~\ree{ga} gives
$$
\ga(\Th_{n,3})
\le
\sum_{d\ge0} (d+1) 2^{3n-2-d}
=
2^{3n-2}\cdot\frac{1}{(1-1/2)^2}
=
2^{3n}.
$$
So we must also have $\ga(\Th_{n,3})\sim 2^{3n}$ since $\ga$ is an
upper bound.

For $\Th_{2,t}$ note that $k$, the number of edges not in a spanning
tree, satisfies $k=t-1$.  We also have $p=t+2$ and $q=2t$.
Using~\ree{ld}, we get
$$
l(d,t)\stackrel{\rm def}{=}l_d(\Th_{2,t})=
\left|L_d\left(x_1^2 x_2^2\cdots x_{t-1}^2\right)\right|
$$
which is the coefficient of $y^d$ in the expansion of
the generating function $(1+y+y^2)^{t-1}$.  From this, it follows that
the $l(d,t)$ satisfy the recursion
\beq
\label{rr}
l(d,t+1)=l(d,t)+l(d-1,t)+l(d-2,t).
\eeq
Let $\ga_t=\ga(\Th_{2,t})$.  So multiplying~\ree{rr} by $2^{t+2-d}$
and summing over $0\le d\le t+2$, we can use~\ree{ga} to get the
following equation, with the three expressions in brackets coming from
the three terms of the recursion (respectively):
\beq
\label{gat}
\ga_{t+1}=
\left[2\ga_t + l(t+2,t) \right]+\left[\ga_t\right]
+\left[\frac{1}{2}\ga_t - \frac{1}{2} l(t+1,t)\right]
>\frac{7}{2}\ga_t - \frac{1}{4} \ga_t = \frac{13}{4}\ga_t
\eeq
where the inequality follows by noting $4 l(t-1,t)$ is a summand in
$\ga_t$ and that, as provable from generating function, the sequence 
$(l(d,t))_{0\le d\le 2t-2}$ is symmetric and unimodal with maximum at $l(t-1,t)$.

Finally, combining the estimates in~\ree{be} and~\ree{gat}, we see that
for any $0<\ep<1/4$,
$$
\be(\Th_{2,t})=(t+1)^2 3^t \exp\frac{-2t}{6t+6}=o((13/4-\ep)^t)=
o(\ga(\Th_{2,t}))
$$
as desired.
\qqed

Now for $n\ge4$, let $\Phi_{n,t}$ be a phi graph derived by pasting
together $t$ 
cycles of length $n$ in such a way that each cycle only intersects the
cycle just preceding and the cycle just following it (if any).  Note
that $\Phi_{n,t}$ is actually a graph 
family  since one can get a number of graphs with these
specifications by pasting along different edges.  But they all have a
uniform description of their NBC bases and degree sequences, so the
bounds under consideration will apply to any graph of the family.

\bth
As $n\ra\infty$ we have
$$
\ga(\Phi_{n,2})\sim\al(\Phi_{n,2}).
$$
As $t\ra\infty$ we have
$$
\ga(\Phi_{4,t})=o(\be(\Phi_{4,t})).
$$
\eth
\pf\
The proof for $\Phi_{n,2}$ is completely analogous to the proof given
for $\Th_{n,3}$, so we leave it to the reader.

Now considering $P_2\times P_{t+1}$ or any other member of $\Phi_{4,t}$,
we see that $p=2t+2$, $q=3t+1$, and $k=t$.
Using the bound~\ree{ch} and the Binomial Theorem
in~\ree{ga} yields
$$
\ga(\Phi_{4,t})\le\sum_{d=0}^{2t+1} {d+t-1\choose t-1} 2^{2t+1-d}
\le 2^{2t+1}\sum_{d=0}^{\infty} {d+t-1\choose t-1} 2^{-d}
= 2^{2t+1}\frac{1}{(1-1/2)^t}=2\cdot 8^t.
$$
Now~\ree{be} gives
$$
\be(\Phi_{4,t})\sim a\cdot b^t,\quad b\approx 14.5682
$$
finishing the proof of the theorem.\qqed

\section{Comments and Open Problems}
\label{cop}

\subsection{Arbitrary fields}

We will now indicate how to generalize our construction to an
arbitrary field.  We first need to review what Brown's hsop looks like
over a field $F$.  Fix an orientation of $E(G)$.  Also, for each $e_j\in
T(G)$, orient all the edges of $D_j$ in one of the two possible
directions.  Now define signs
$$
\ep_{i,j}=
\case{1}{if the orientation of $e_i$ in $G$ is the same as in $D_j$,}
{-1}{if these orientations are opposite.}
$$
We have corresponding polynomials
\beq
\label{thF}
\th_j=\sum_{e_i\in D_j} \ep_{i,j} x_i.
\eeq
\bth[\cite{bro:cpo}]
If $G$ is a connected graph then 
the set of polynomials defined by~\ree{thF} for $e\in T(G)$ is an hsop for
$F(G)$.\qqed 
\eth

Solving for $x_j$ in the equation for $\th_j$ and plugging into the
monomials $\bx^\Cb$, $C\in\cC(G)$, gives the generators $p_\Cb$ for an ideal
$J(G)$ such that
$$
R(G)\iso \frac{F[x_1,\ldots,x_k]}{J(G)}.
$$
Note that the monomial $m_\Cb$ that was chosen from the expansion of
$p_\Cb$ in the case $F=\bbZ_2$ will also appear with coefficient
$\pm1$ for any field.  So the proof of Theorem~\ref{deg2} will go
through as before as long as the generators of $J(G)$, $J(G\sm e_q)$,
and $J(G/e_q)$ can be related in the correct way.

An orientation of $G$ induces orientations of $G\sm e_q$ and $G/e_q$
merely by keeping each $e_i$, $i<q$, oriented the same way in all
three graphs.  Under the assumptions of Theorem~\ref{deg2} we showed
that $D_j(G\sm e_q)=D_j(G)-e_k$ for $j>k$.  So we can orient 
$D_j(G\sm e_q)$ the same way as $D(G)$ in this case.  We also have
$D_k(G\sm e_q) = \{e_k\}$, so it does not matter which way we orient
$e_k$ in this cut set as $x_k$ is being set to zero in the quotient.
Thus we get, as we did in the $\bbZ_2$ proof, that the generators for
$J(G\sm e_q)$ are gotten from those for $J(G)$ by setting $x_k=0$.
Similar considerations show that we can define orientations on the cut
sets of $G/e_q$ so that the equalities we had before still hold.  So
Theorem~\ref{deg2} holds, and hence so do all the rest of the results of the
previous sections, over any field.

\subsection{Arbitrary graphs}

We conjecture that any graph $G$, with its edge set suitably ordered,
has an NBC basis.
\bcon
\label{C}
Let $G$ be any graph.  Then there is a standard ordering of
$E(G)$ such that $L(G)$ is a basis for $R(G)$. 
\econ

We will now outline a possible line of attack on Conjecture~\ref{C}.
Even though we have not been able to push it through, it is possible
that some of these ideas will be useful in finally proving or
disproving this conjecture.
Recall that it suffices to find a proof when $G$ is a block.
But any block other than $K_2$ (the complete graph on 2 vertices)  
has a nice recursive structure in that it can be built
from a cycle by adding a sequence of paths called {\it ears\/}.  This
result is due to Whitney~\cite{whi:cgc}.  Proofs can also be found in
the books of Diestel~\cite[Proposition 3.1.2]{die:gt} and
West~\cite[Theorem 4.2.8]{wes:igt}.
\bth[Ear Decomposition Theorem]
Suppose $G\neq K_2$.  Then
$G$ is a block if and only if there is a sequence
$$
G_0, G_1,\ldots, G_l=G
$$
such that $G_0$ is a cycle and $G_{i+1}$ is obtained  by
taking a nontrivial path and identifying its two endpoints with two
distinct vertices of $G_i$. \qqed
\eth

Note that the graph $G_1$ in the ear decomposition sequence is a theta
graph.  So one might try to prove Conjecture~\ref{C} by induction on
$l$, the number of paths added.  (Actually, one also needs to induct on the
number of edges since one contracts an edge and not a whole path.)  In
fact, the induction step goes 
through in much the same way as our proof for theta graphs as long as
the path added has length at least two.  The difficulty comes if the
path is a single edge.  In that case, it is still easy to relate the
circuits of $G\sm e_q$, where $e_q$ is the newly added edge, to those
of $G$.  But the situation is much more complicated in $G/e_q$, which
may not even be a block.  So a more delicate analysis is needed.  
Unfortunately, there are graphs (such as the complete graphs) where every ear
decomposition requires the addition of a single edge at some stage.

\subsection{Not quite arbitrary matroids}

As a last point, the reader may have noticed that all of the graphical
definitions we used to define NBC bases make sense for the broken
cirucit complex of an arbitrary
matroid.  So a natural question is whether our construction goes through in that
level of generality.  Brown, Colbourn, and Wagner~\cite{bcw:cmr} have
a way of producing an hsop for any representable matroid.  (Actually,
their construction is of an hsop for the independence complex of the
matroid.  But this will also give an hsop for the broken circuit
complex since it is a subcomplex of the independence complex having the
same rank.)  So this
would be the natural class of matroids in which to look for NBC bases.

\subsection{Gr\"obner bases}

We note that, in general, the monomials used to generate $U(G)$ are not the leading
terms of a  Gr\"obner basis for the ideal $J(G)$.
As an example of this, one can take a theta graph
consisting of three paths of length two in the theta labeling as
described in Section~\ref{tf}.  Then by choosing a suitable
orientation for $G$ and its cocircuits, $J(G)$ will have generators
\beq
\label{gen}
\{p_\Cb\}=\{x_1(x_1+x_2)^2,\ x_2(x_1+x_2)^2,\ x_1x_2^2\}
\eeq
from which we pick monomials
\beq
\label{mon}
\{m_\Cb\}=\{x_1^3,\ x_2^3,\ x_1x_2^2\}
\eeq
for the NBC basis.

Suppose, towards a contradiction, that there is a term ordering
giving~\ree{mon} as the set of leading terms of a Gr\"obner basis.  Then in
that term ordering we either have $x_1<x_2$ or $x_1>x_2$.  Suppose the
former is true.  Then $x_1^3$ is the smallest (monic) polynomial which
is homogeneous of degree three.  Also, the generators of $J(G)$ are
homogeneous.  So if $x_1^3$ were a leading term of a polynomial in
$J(G)$ then, in fact, $x_1^3\in J(G)$.  But it is easy to check that
$x_1^3\not\in J(G)$ since it is not a linear combination of the
polynomials in~\ree{gen}.  (It suffices to consider linear
combinations by homogeneity.)  One gets a similar contradiction using
$x_2^3$ if one assumes that $x_1>x_2$.  So no such Gr\"obner basis can exist.

\medskip

{\it Acknowledgement.\/}  We would like to thank David Forge for
stimulating discussions.

\bigskip
\bibliographystyle{plain}
\bibliography{ref}

\end{document}